\documentclass[letter,12pt]{article}
\usepackage{amsfonts,amsmath,amssymb}
\usepackage{mathptmx}       % selects Times Roman as basic font
\usepackage{helvet}         % selects Helvetica as sans-serif font
\usepackage{courier}        % selects Courier as typewriter font
\usepackage{type1cm}        % activate if the above 3 fonts are
\usepackage{makeidx}         % allows index generation
   \usepackage[dvips]{graphicx}
 
   \DeclareGraphicsExtensions{.eps}

\DeclareMathOperator{\rank}{rank}
\newcommand{\tx}{\tilde{\mathbf x}} 
\newcommand{\ty}{\tilde{\mathbf y}} 
 
\newcommand{\tF}{\tilde{\Phi}} 
\makeindex             % used for the subject index
\begin{document}

\title{Correcting Errors in Linear Measurements and  Compressed Sensing of Multiple Sources.}
% Use \titlerunning{Short Title} for an abbreviated version of
% your contribution title if the original one is too long
\author{Alexander Petukhov and Inna Kozlov}
% Use \authorrunning{Short Title} for an abbreviated version of
% your contribution title if the original one is too long
%\institute{Alexander Petukhov \at University of Georgia, Athens, GA 30602, \email{petukhov@math.uga.edu}
%\and Inna Kozlov \at Algosoft Tech, Bogart, GA \email{kozlovinna@yahoo.com}}
%
% Use the package "url.sty" to avoid
% problems with special characters
% used in your e-mail or web address
%
\maketitle

\abstract{We present an algorithm for finding sparse solutions of the system of linear equations
$\Phi\mathbf{x}=\mathbf{y}$ with rectangular matrices $\Phi$ of size $n\times N$, where $n<N$, when measurement vector $\mathbf{y}$ is corrupted by a sparse vector of errors $\mathbf e$. 
\par
We call our algorithm the $\ell^1$-greedy-generous (LGGA) since it combines both greedy and generous strategies in decoding.
\par
Main advantage of LGGA over traditional error correcting methods consists in its ability to work efficiently directly on linear data measurements. It uses the natural residual redundancy of the measurements  and does not require any additional redundant channel encoding.
\par
We show how to use this algorithm for encoding-decoding  multichannel sources. This algorithm has a significant advantage over existing straightforward decoders when the encoded sources have different density/sparsity of the information content. That nice property can be used for very efficient blockwise encoding of the sets of data with a non-uniform distribution of the information.  The images are the most typical example of such sources.  
\par
The important feature of  LGGA  is its separation from the encoder. The decoder does not need any additional side information from the encoder except for linear measurements and the knowledge that those measurements created as a linear combination of different sources.
}

\section{Introduction}

The problem of finding sparse solutions of a system of linear
equations 
\begin{equation}
\label{sle}
\Phi\mathbf x=\mathbf y,\quad \mathbf x\in\mathbb R^N, \,\mathbf y\in\mathbb R^n,\,N>n,
\end{equation}
is interesting in many Information Theory related contexts. It is tractable as reconstruction of a~sparse data 
vector $\mathbf x$ compressed with the linear operator $\Phi$. While from the point of view of classical linear algebra,
system (\ref{sle}) may not have a unique solution, the regularization of the problem in the form of the solution sparsity allows to guarantee the uniqueness (or the uniqueness almost for sure) in many practically important cases. 
\par

 Within this paper, we say that a vector $\mathbf a\in \mathbb R^m$ is {\it sparse} if $$k:=\#\{i\mid a_i\ne 0,\, 1\le i\le m\}<m,$$ i.e., it has some zero entries. This number is called also the Hamming weight of the vector $\mathbf x$. For a randomly selected  matrix $\Phi$ and a vector $\mathbf y$, if the~sparse solution  of system (\ref{sle}) exists,  
 it is unique almost for sure. In particular, a sparse, i.e., having at most $n-1$ non-zero entries, vector $\mathbf x$ theoretically can be restored from its measurements $\mathbf y$ almost for sure.
\par
Unfortunately, a straightforward search for sparse solutions to (\ref{sle}) is an NP-hard problem (\cite{Nt}). The NP-hardness of the problem does not contradict the existence of the algorithms for finding sparse solutions in some quite typical cases.
\par
The recent Compressed Sensing (Compressive Sampling) studies gave a big push for  the development of the theory and  affordable algorithms for finding sparse solutions. 
It turned out that for a reasonable (not very large) value of the~sparsity the vector $\mathbf x$ can be recovered precisely, using
Linear Programming Algorithm (LPA) for finding the solution to (\ref{sle}) with the minimum of $\ell^1$-norm (\cite{CT}, \cite{D},
\cite{RV}). While, in practice, the number $k$ characterizing the sparsity of the vectors $\mathbf x$ which can be recovered with $\ell^1$-minimization is far from the magic number $n-1$, this case is a well-studied and reliable tool for solving systems (\ref{sle}) for many applied problems. 
\par
Orthogonal Greedy Algorithm (OGA) is a strong competitor of LPA. If it is implemented appropriately 
(\cite{DTDS}, \cite{P}),  it outperforms LPA in both the computational complexity and in the ability to recover sparse representations with the~greater sparsity number $k$. In some papers (e.g., \cite{DTDS}), this modification of OGA is called Stagewise Orthogonal Matching Pursuit (StOMP = StOGA).
\par
Very resent success in the efficient recovery of sparse representation is due to paper \cite{KMSSZ}.
The authors suggested to use a special (band diagonal) type of sensing matrices in combination with 
a belief propagation-based algorithm.
While the algorithm \cite{KMSSZ} gives a very powerful tool  in linear methods of data
compression, it is oriented on a special form of measuring matrix which seems to be theoretically less efficient from the point of view of the error correcting capability.
\par
We do not discuss hear the  advantages and disadvantages of different decoding algorithms in detail. For the goals of this paper, we are interested in algorithms supporting the opportunity  to assign non-equal significance of the solution entries. Such algorithms use the idea of reweighted / greedy optimization (see \cite{CWB}, \cite{CY}, \cite{DDFG}, \cite{FL}, \cite{KP}, \cite{KP1}). In what follows, we use our $\ell^1$-greedy algorithm (LGA, see Algorithm A below) as a starting point for the algorithms considered in this paper. LGA has 2 advantages over competitors. First of all, it was shown numerically in \cite{KP} that LGA has the highest capacity of the recovery of sparse/compressible data encoded  with Gaussian matrices. Second, what is more important, it is easily adaptable 
to the needs of this paper. 
\par
The absence of theoretical  justification and relatively high computational complexity are main disadvantages of LGA.  However, it should be emphasized that the main competitors of LGA also
do not have appropriate theoretical justification. As for the computational complexity, the fast version of LGA was 
developed in \cite{KP1}. Its computational complexity is about the complexity of the regular 
$\ell^1$-minimization, whereas the reconstruction capacity is very close to regular LGA and significantly higher than other, even more computationally extensive, algorithms.
\par
Above and in what follows, saying that the algorithm has "the higher recovery capability", we mean that for the same $n$ and $N$ the algorithm is able to recover vectors $\mathbf x$ with a larger number $k$ of non-zero entries. Of course, the relative capability of algorithms depends on the statistical model of non-zero entries. We use only the Gaussian model since a Gaussian random value has the highest information content among  random values with the same variance.  The reweighted algorithms are known as  inefficient for  vectors $\mathbf x$ with th Bernoulli distribution (equally probable $\pm 1$ for non-zero components). However, we believe that CS for Bernoulli input has the less applied value since Bernoulli distribution has a very low information content and the linear methods of compression are extremely inefficient for them.  
At the same time, for distributions having the sparsity (the decay rate of non-zero entries) higher than for the Gaussian distribution the capability of greedy-based algorithms is higher than for the Gaussian distribution.
\par
The state-of-the-art approach to data encoding is based on the concept (due to C. Shannon) that the optimal 
data encoding for transmission through a lossy channel can be splitted into 2 independent stages 
which are source and channel encoding. 
\par
The first stage is source encoding or compression. Usually, compression is a non-linear operation. In the ideal case of the optimal compression, its output is a bitstream consisting of absolutely random bits with equally probable "0"s and "1"s.  The compression reduces the {\it natural} redundancy of the encoded data.
\par
The second stage is channel encoding. It  plays the opposite role, introducing the {\it artificial} redundancy. Usually (but not necessary), this is a linear operator on a Galois (say binary) field. The redundancy introduced here is different from the redundancy removed on the compression stage. Its model is completely known to the receiver of the information (the decoder) and in the case of moderate corruption this model can be used for the perfect recovery 
of the transmitted data. 
\par
While the natural redundancy of the source also can be used for error correction, it is not so reliable for this purpose as channel encoding. Nevertheless, in some deeply studied problems a kind of "error correction" based on the natural redundancy is possible. Among those applied problems we just mention various data denoising methods and digital film restoration.  The natural redundancy arises when  a data model does not allow the  digital data representation to take arbitrary values. For the denoising the belonging of the data to some class of smoothness serves for protection against corruption. Whereas impossibility of big changes between frames following with the rate 20 -- 30 frames per second plays the same role for moving images recovery. 
In both cases, the redundancy model is known only approximately. 
\par
One more problem indirectly conneted to the error correction is image upsampling. The inpainting of missing information is tractable as correcting "errors"  lost in a  channel with erasures.
\par
In CS community, applications to error correcting codes were noticed practically simultaneously with the mainstream compression issue (cf., \cite{CR}, \cite{CRTV}, \cite{RV}). Repeating commonly accepted channel encoding strategy, 
this approach consists in introducing the redundancy by  multiplying the data vector $\mathbf y$  (say, the "compressed" output in (\ref{sle})) by a matrix $B$ of the dimension $m\times n$, where $m>n$. 
 We assume that $\rank B=n$.
Then the output vector
$$
\mathbf z:=B\mathbf y=B\Phi\mathbf x
$$
is protected, at least theoretically, from the corruption of up to $m-n-1$ its entries.  To understand the mechanism of this protection, introduce the matrix $C$ of dimension $m\times(m-n)$ whose columns constitute a basis of the orthogonal complement of the space spanned by the columns of $B$. Assuming that $\mathbf e$ is a sparse vector of errors,  compute the vector  $\mathbf s:=C^T(\mathbf z+\mathbf e)$,   known in the Coding Theory as a syndrome. The syndrome can be measured by the receiver from corrupted information. At the same time, $C^T\mathbf z=C^TB\mathbf y=\mathbf 0$, the corruption vector $\mathbf e$ can be found as a sparse solution of the system
$$
C^T\mathbf e=\mathbf s.
$$
\par
The scheme considered in the previous paragraph is in the intersection of mainstreams of CS and Coding Theory based on separation of the source and channel encoding. 
\par
In this paper, we discuss a different data protection scheme which uses the residual space not occupied by the data for error correction purposes.  In the  literature on Coding Theory, a data transform providing simultaneous compression and protection from transmission errors is called  Joint Source-Channel Coding.  As we will see below, the plain linear measurements have
this property provided that the encoded data are sparse enough in some representation system and the vector of corruption is also sparse. Therefore, in this case, the compressed data are restorable without any additional channel encoding. To our knowledge,  the first time, this idea was formulated
in \cite{WYGSM}.
\par
In this paper, we do not consider any methods for fighting the corruption with the noise in entries when the level of the corruption of the output vector $\mathbf y$ is relatively low
 but the corruption takes place at each entry. However, we will discuss the stability of our algorithm 
with respect to such corruption.
\par
We emphasize that the encoding model accepted in this paper is analog encoding, i.e., we do not mean to
apply analog-to-digital transform to the measurements.  The entries are stored (transmitted) in the form
of their magnitude (not in bitwise representation).
 For this reason we concider the sparse model of corruption when errors  are introduced directly into entries of the vector $\mathbf y$ (not in bits).
\par
The structure of the  paper is as follows. In Section \ref{SecProblemSetting}, we descuss how  corrupted data can be recovered with the regular LGA algorithm. In Section \ref{SecVarErrorRate}, we show that the "generous" (LGGA) modification of LGA provides much higher error correcting capability provided that the error rate is known at least approximately.   In Section \ref{SecMultiSource}, using the well known theoretical equivalency 
of  a lossy channel and multisource transmission, we discuss advantages of joint encoding of a few sources of data.  In Section \ref{SecAdaptAlg}, we design the adaptive $\ell^1$-greedy-generous algorithm (ALGGA) which does not require preliminary knowledge
of either error rate or the relative density of information in multiple sources.
\par
While all results of the paper are based exceptionally on  numerical experiments, the stability of those
results allow to be optimistic about possible applications of algorithms designed on the greedy-generous  principles introduced below.
\par
In all our numerical experiments we used random matrices $\Phi$ of size $128\times 256$. This size seems to be the most popular in academic papers. The selection of other parameters which can be found
in the text below  is associated with that size and its interplay with LGA. They make sense in association with
the size of input data. While we made experiments with different sizes of $\Phi$, we decided to do not include those 
results because they just confirm the effects and efficiency of the algorithms considered in the paper and do not bring new effects deserving extra journal space.
\section{Problem setting and approaches to solving}
\par
\label{SecProblemSetting}
Let $\mathbf x\in\mathbb R^N$ be a sparse vector with $k$ non-zero entries. 
The vector  $\mathbf x$ is encoded with a linear transform $\Phi\mathbf x=:\mathbf y\in\mathbb R^n$. The vector $\mathbf y$ is corrupted with a vector $\mathbf e\in\mathbb R^n$ with $r<n$ non-zero entries. 
The decoder receives the corrupted measurements
$$
\tilde{\mathbf y}=\Phi\mathbf x+\mathbf e.
$$
The matrix $\Phi$ is also available for the decoder.
It is required to find the data vector $\mathbf x$. Of course, when $\mathbf x$ is found, the error vector also can be computed.
Easy to see that the problem can be rewritten in the form 
\begin{equation}
\label{extsys}
\ty=\tF\tx, \, \text{ were } \tF:=(\Phi\,\, I_n)\in\mathbb R^{n\times(N+n)}, \,\tx=
\left(
\begin{array}{l} \mathbf x \\ \mathbf e\end{array}
\right),
\end{equation}
$I_n$ is the identity $n\times n$ matrix. Assuming that the columns of the matrix $\Phi$ are almost orthogonal to columns of $I_n$ and taking into account that both $\mathbf x$ and $\mathbf e$ are sparse vectors and 
$k+r<n$, we arrive at a standard problem of finding a sparse solution of an underdetermined system of linear equations with the matrix $\tF$. We emphasize that although the columns of $I_n$ and 
$\Phi$ span the same space $\mathbb R^n$ and do not introduce separation between data and error spaces, the spans of the columns corresponding to the sparse vectors $\mathbf x$ and $\mathbf e$ are approximately orthogonal, so the errors can be separated from the data.
\par
Thus, the error correction problem is reduced to the problem which can be solved with any standard CS decoder. 
The algorithm we are goint to use for decoding is based on iterative minimization of the functional
\begin{equation}
\label{weighted}
\|\mathbf x\|_{w_{j,1}}:=\sum w_{j,i}|x_i|,
\end{equation}
which is the weighted $\ell^1$-norm. The initial idea to use the reweighted $\ell^1$-minimization for sparse solutions of underdetermined system is 
due to \cite{CWB}. We will use the following algorithm with the higher recovery capability.
\par
{\bf Algorithm A}
({\it the $\ell^1$-Greedy Algorithm}  (LGA)\cite{KP})
\par
\begin{enumerate}
\item {\it Initialization}.
\begin{enumerate}
\item
 Take take $w_{0,i}\equiv 1$, $i=1,\dots,N$.
\item
 Find a solution  $\mathbf x_0$ to (\ref{sle}) providing the minimum for (\ref{weighted}).
\item
 Set $M_0:=\max\{x_{0,i}\}$ and $j:=1$.
\end{enumerate}
\item
{\it General Step}.
\begin{enumerate}
\item Set $M_j:=\alpha M_{j-1}$ and the weights
$$
w_{j,i}:=\left\{\begin{array}{ll}\epsilon, &|x_{j-1,i}|>M_j,\\ 1, &|x_{j-1,i}|\le M_j.\end{array}\right.
$$
\item
 Find a solution  $\mathbf x_j$ to (\ref{sle}) providing the minimum for (\ref{weighted}).
\item
 Set $j:=j+1$.
\item
If Stopping Criterion is not satisfied, GOTO 2 (a).
\end{enumerate}
\end{enumerate}
\par

\par
LGA has a minor difference with the standard reweighted $\ell^1$-minimization from \cite{CWB}. This difference consists in dynamic updating the weighting function.  Nevertheless, on the Gaussian input, LGA outperforms  both the regular and the reweighted $\ell^1$-minimization significantly. 
\par
The results of straightforward application of Algorithm A to the extended inputs $\ty$ and $\tF$ from (\ref{extsys}) will be shown below.
Before we present numerical results, we describe our settings applicable to all numerical experiments in this paper. We generate vectors $\mathbf x$ and $\mathbf e$ with correspondingly $k$ and $r$ non-zero entries with normal distribution with parameters (0,1).
The entries of the matrix $\Phi$ also has taken from the standard normal distribution.
\par
 We run  200 independent trials with Algorithm A or its modifications. In each trial, the result is obtained after 30 iterations of the algorithm with parameters $\alpha=0.85$, $\epsilon=0.001$ or
if the club of large coefficients, whose magnitude exceeds $M_j$, reaches the cardinality $n$. 

\par
\begin{figure}[h]
\centering
\includegraphics[width=3.7in]{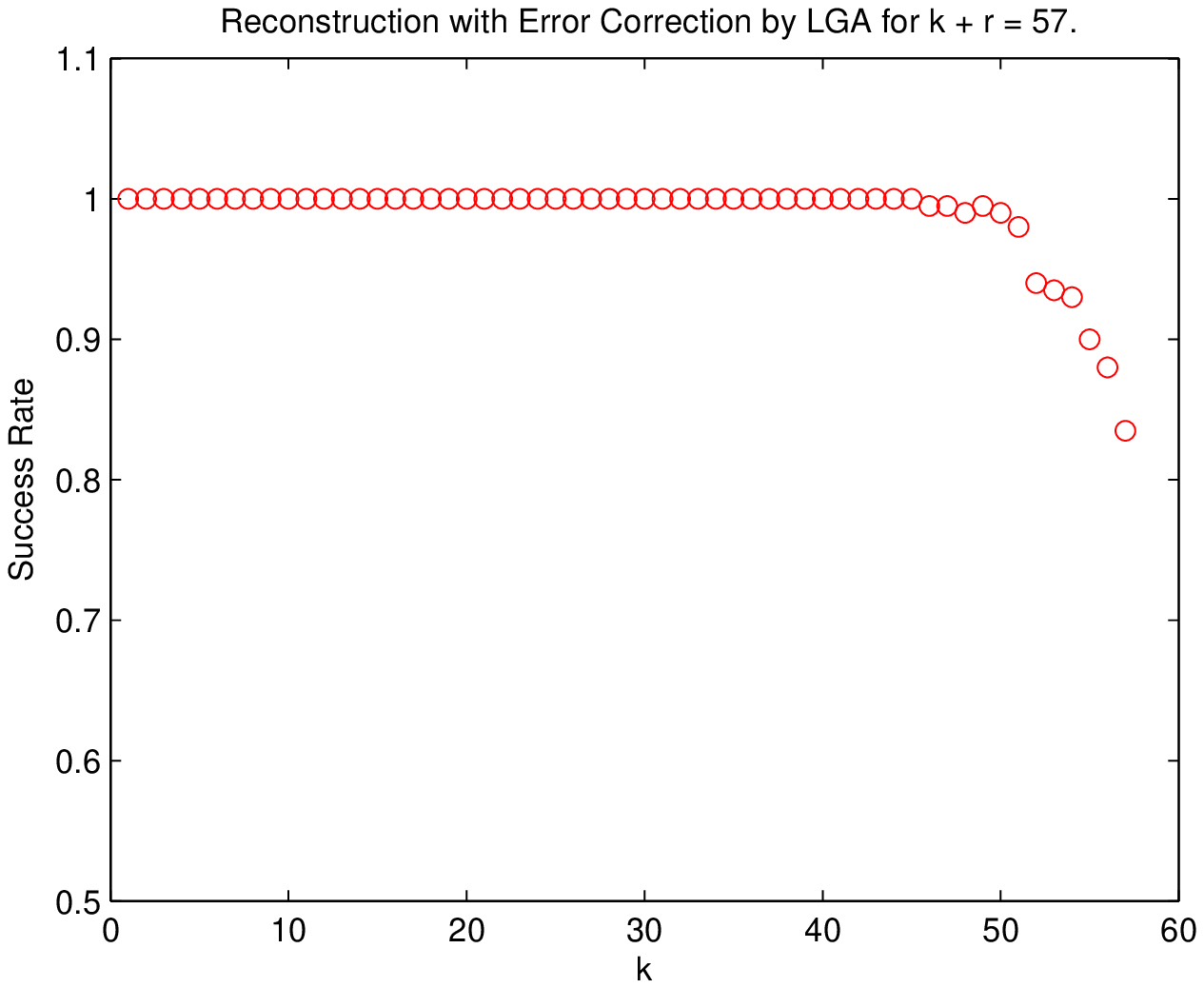}
\caption{Error correction  with Algorithm A for $\Phi\in\mathbb R^{128\times 256}$} when the~total number of non-zero entries and errors is fixed (=57).
\label{fig_KplusRequal57}
\end{figure}

On Figure \ref{fig_KplusRequal57}, we present the graph of the relative success rate of the recovery for $k$-sparse data by Algorithm A when $k+r=57$. The number 57 corresponds to the LGA reconstruction success rate about 0.82 for sensing 57-sparse vectors with random $128\times384$ matrices. The horizontal axis reflects the values of $k$ in the range $1\div 57$, whereas the number of errors for each point of the graph is computed as $r=57-k$. 

\par
Since we have constant the number $57=k+r$ of "information" entries,
in general, the almost constant behavior of the graph on Fig. \ref{fig_KplusRequal57} is  predictable, except for the better performance of the recovery algorithm when we have less non-zero data entries in $\mathbf x$ and a greater number of errors. Such behavior can be explained by the fact that the errors are associated with the identity
part of the matrix $\tF$, which has orthogonal columns, whereas the columns of the matrix $\Phi$ are only approximately orthogonal.  Therefore errors can be found and corrected easier than the data entries. Of course, if the entire extended matrix $\tF$ were constructed with the Gaussian distribution, the graph would be parallel to the $k$-axis.
\par
Thus, Algorithm A works for  the error correction even better than it was expected in advance.
\par
\begin{figure}[h]
\centering
\includegraphics[width=3.7in]{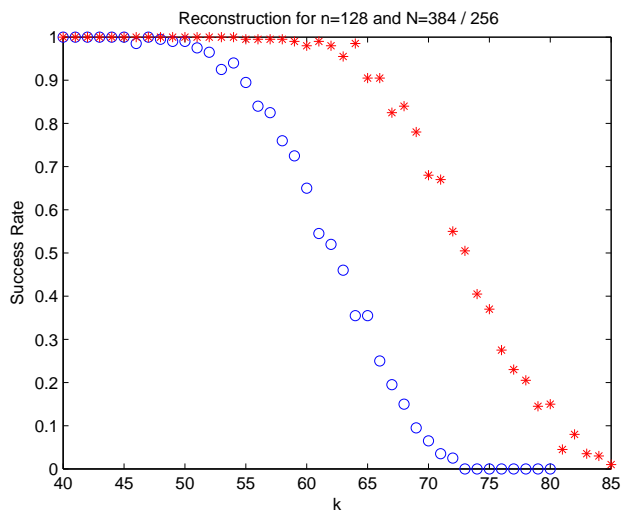}
\caption{Algorithm A (error free) for  $\Phi\in\mathbb R^{128\times 256}$(right)  and $\Phi\in\mathbb R^{128\times 384}$(left)}.
\label{fig_256and384}
\end{figure}

However, there is one very discouraging drawback of such error correction. 
If we run it on a 57-sparse vector and its CS measurements are error free, then the Algorithm A with the extended matrix $\tF\in\mathbb R^{128\times 384}$ will recover it with the probability about 0.82, whereas the same LGA algorithm applied to the matrix $\Phi\in\mathbb R^{128\times256}$ is able to recover vectors of the  same sparsity with probability very close to 1. 
The probability 0.82 on $\Phi$ is reached by Algorithm A at the sparsity $k=68$.
Thus, we loose about 20\% of the efficiency just because of the suspicion that the data could be corrupted.  The curves of the reconstruction success rates for $N=256$ and $N=384$ when $n=128$ which demonstrate those losses are presented on Fig. \ref{fig_256and384}.
\par
\par
From the point of view of the multidimensional geometry everything looks reasonable. We use $k$-dimensional subspaces of $\mathbb R^n$ which are all possible spans of $k$-subsets of a set of $N$ vectors from $\mathbb R^n$.  If we increase that set by $n$ vectors, then the number of $k$-planes  drastically increases from $N\choose k$ to ${N+n}\choose k$. Therefore the $k$-planes become harder distinguishable. 
\par
From the point of view of classical coding theory those losses also are explainable. The error  resilience requires some redundancy in the representation. The higher projected error rate the bigger redundancy has to be introduced in advance. If we are lucky and a channel is error free, the redundancy is already spent (actually, in this case, it is wasted) on the encoding stage. Any improvement of the recovery due to the absence of the errors is impossible.
\par
There is a fundamental difference of our situation. We do not introduce any artificial redundancy. Therefore, for the lower rate of errors we may hope for the improved performance of the decoder on the encoded data recovery. 
The next section explains how this intuitive advantage can be transformed into actual benefits of the decoder.
%%%%%%%%%%%%%%%%%%%%%%%%%%%%%%%%%%%%%
\section{Algorithm for Variable Error Rate}
\label{SecVarErrorRate}
%%%%%%%%%%%%%%%%%%%
Let us come back to Fig. \ref{fig_256and384}
to specify what kind of error correcting algorithm we want to design. 
If we apply Algorithm A to the potentially corrupted vector $\ty$ obtained as a measurement with the extended matrix $\tF$, then the success rate will be described by the curve corresponding to the matrix $128\times(256+128)$ (i.e., $N=384$). At the same time, provided that we know  that there are no errors in the channel, using this piece of information, we just switch Algorithm A to the matrix $\Phi$ of size $128\times 256$ ($N=256$) and get much higher success rates. Thus, Algorithm A is trivially "adaptable" in this way to the case $r=0$.
\par
For $r=1$ the situation changes drastically. If we apply Algorithm A to $\tF$,
this curve is just a translation of the curve for $N=384$ by 1 to the left what is close to the error free case.  While applying Algorithm A to the matrix
$\Phi$, we do not get sparse solutions at all because the columns of the identity matrix do not have sparse representations in the dictionary of columns of the matrix $\Phi$. In this case, even if we have side information that only one entry of  the data are corrupted, Algorithm A applied to $\ty$ and $\Phi$ becomes  forceless,  whereas applying it to $\tF$ we lose in efficiency significantly. We wonder whether we can modify Algorithm A in such way that the success rate curve becomes the translate of the curve for $N=256$ by $r=1$ to the left rather than the curve for $N=384$. 
\par
This goal can be reached with Algorithm A if we solve 128 problems with matrices $\Phi_i\in \mathbb R^{128\times257}$, $i=1, 2,\dots,128$, obtained as an extension of $\Phi$ with the $i$th column of the identity matrix and select the sparsest of 128 solutions. Such algorithm is computable but
very computationally extensive and, because of this complexity, cannot be extended to $r$ even slightly greater than 1.
\par
Thus, the main question of this sections is how the information about the density of errors can be incorporated in Algorithm A for achieving the maximum possible success rate with minor or no change of computational complexity. 
\par
Let us discuss the potential limits of the desired algorithm. First of all, when 
\begin{equation}
\label{balance}
\frac{k}{N}\approx\frac{r}{n},
\end{equation}
i.e., the density of  actual errors in $\mathbf e$ is equal to the density of non-zero entries of $\mathbf x$,  and non-zero entries of $\mathbf x$ and $\mathbf e$ have same distributions of their magnitudes, we cannot hope for increasing the efficiency, staying within the framework of  LGA. 
Indeed, under those assumptions, we can think about the vector $\tx$ as about a sparse vector with a randomly and independently selected index set of $k+r$ non-zero components among $N+n$ entries encoded as $\ty$ with the matrix $\tF$. Definitely, the point $k+r$ will rather correspond to the point of the curve for $N+n$ (= 384 above) than $N$ (= 256 above). Thus, we cannot expect any improvement  in the recovery of the vector $\tx\in\mathbb R^{N+n}$ with respect to the success rate curve for $N+n$. 
Of course, as we noticed above the result can be slightly better just because of the structure of $\tF$ involving the identity matrix.
\par
At the same time, we have a right to expect the improvement when the equality of the proportions (\ref{balance}) is violated. For the case $k/N\gg r/n$ we may expect the success rate for the recovery $k+r$ non-zero entries close to the success curve corresponding to $n\times N$ matrix, whereas for $k/N\ll r/n$, we want to see the recovery rate close to the recovery rate for the identity matrix $I_n$. 
\par
At the first glance, the last request may look too challenging because the identity matrix encoder $\mathbf y:=I_n\mathbf x$ allows the trivial "recovery" $\mathbf x:=\mathbf y$. So the decoder must  be extremely efficient. However, it is quite realistic. Indeed, when we encode $\mathbf x$ whith only one non-zero entry (i.e., $k=1$) we can recover it with a big probability even if $\mathbf e$ has up to $n-2$ non-zero entries. The direct exhaustive search for the 1-sparse solution in this case has the qubic (actually, $\sim n^2N$) complexity. This complexity as well as the required precision are high but still realistic. Therefore the high eficiency of the projected algorithm for $r$ close  to $n$ would not be considered as a miracle.
\par
Let us formulate our goals in the developing a new recovering algorithm with error correcting capability. The algorithm has to provide the property which has meaning of the 3-point interpolation. It has to work at least as Algorithm A on compressing matrices of size $n\times(N+n)$ when (\ref{balance}) takes place. Its efficiency has to approach the cases of  compressing matrices of size $n\times N$ and $n\times n$ when $r\to 0$ and $k\to 0$ correspondingly.
\par
The idea of such algorithm relies on  the reweighting strategy and some knowledge about an approximate error rate. 
\par
Algorithm A is based on  the greedy  strategy. Its main idea is to provide the entries with the bigger expected magnitude with more freedom in the representation of the vector $\ty$. This freedom is provided by the significantly less weight in the (weighted) $\ell^1$-norm given to the entries assumed to be large in the ``optimal''  (sparse) representation of $\ty$. Then the algorithm does not worry too much about their contribution into the norm. So they can be used efficiently for  the partial sum minimizing the residual.
This is a typical  greedy strategy when large coefficients are selected in a separated set and on the next iteration the biggest coefficients in the representation of the residual extend ``the club of large coefficients''.
\par
The greedy approach inspired us on the opposite ``generous'' strategy. If we know that with large probability the channel is (almost) error free, let us give  "generously" a larger weight to the entire block of entries corresponding to the errors or, vice versa, we set a less weight to the entire error block when its density is higher than the density of the data entries.
For instance, if we want to modify LPA algorithm according to the generous principle formulated above, instead of the minimization problem
$$
\|\tx\|_1\to \min, \text{ subject to } \tF\tx=\ty,
$$
we solve the problem
$$
\|\mathbf x\|_1+\lambda \|\mathbf e\|_1\to \min, \text{ subject to } \tF\tx=\ty,
$$
where $\lambda$ depends on the density of errors. 
The $\ell^1$-greedy algorithm from \cite{KP} (and its accelerated version from \cite{KP1}) 
reweights input data of its iterations, according to the output of its previous iterations.
The modification of the problem above does not require any significant changes in Algorithm A itself. We just introduce a different weight for the error components not included into the club of large coefficients. To be more precise we replace the definition of the weight in item 2a of Algorithm~A
with
$$
w_{j,i}:=\left\{\begin{array}{ll}\epsilon, &|x_{j-1,i}|>M_j;\\ 1, &|x_{j-1,i}|\le M_j,\,i\le N;\\  \lambda, &|x_{j-1,i}|\le M_j,\, i>N.\end{array}\right.
$$
We call this modification {\bf Algorithm B} or {\it the $\ell^1$-greedy-generous algorithm} (LGGA).
\par
\begin{figure}[h]
\centering
\includegraphics[width=3.7in]{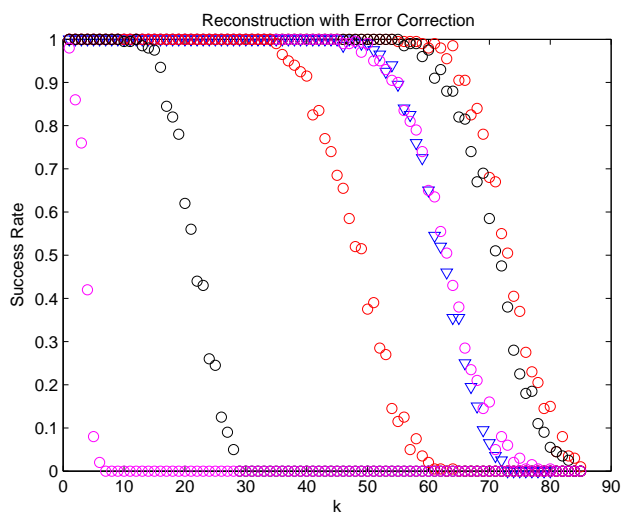}
\caption{Algorithm B for $r=0,1,5,15,45,90$ ("$\circ$"curves from right to left), $\Phi\in\mathbb R^{128\times 256}$ ; \newline
$r=0$, $\Phi\in\mathbb R^{128\times 384}$  ("$\nabla$" curve).}
\label{fig_ErrorCorrectedCurves}
\end{figure}

The efficiency of Algorithm B is illustrated on Fig. \ref{fig_ErrorCorrectedCurves}. We found the curves of the success rate for $r=0,\,1,\,5,\,15,\,45,\, 90$ (the graphs plotted with "$\circ$"s from right to left), applying $\lambda=2.0.\,2.0,\,1.7,\, 1.5, \,0.7, \,0.55$ in Algorithm B. 
%The graph for $r=0$ ($\lambda=2.0$) is drawn with "*"s.

\par
\par
On the same figure, we plot 
the graph corresponding to the error free case for $N=384$ (the graph plotted with ``$\nabla$''s).
We see that that the curve of the success rate of Algorithm A for error free input practically coincides 
with the success rate of Algorithm B for the case $r=5$. So, under similar input conditions, Algorithm B corrects up to 5 errors when Algorithm A admits only error free input.
\par
To discuss the next issue let us move,  in the mind's eye,  each of graphs on Fig.  \ref{fig_ErrorCorrectedCurves} plotted with ``$\circ$''s to the right by the corresponding value $r$. Then we have graphs of the recovery $k+r$ (data and errors) entries. We  see that for $r$ equal to 15 and 45 (and definitely for $r$ between them) the shifted graphs are very close to the graph for the error free case when $N=384$. In this case, generous reweighting does not play any significant role and the results for Algorithm A and B are quite close. For the case $r=1$, the shifted graph is practically coincide with the error free case when $N=256$.  The translate of the curve $r=5$ by 5  brings the graph between error free graphs for $N$ equal to 256 and 384. We note that for the case $r=90$ the shifted graph will be far to the right from the curve $N=256$. This is completely predictable because the identity part of $\tF$ is not compression at all. Therefore, for the matrix $I_n$ separated from $\Phi$, any number of entries are ``recoverable'' with the trivial "algorithm" $\mathbf x:=\mathbf y$. 
\par
If we compare the graphs on Fig. \ref{fig_ErrorCorrectedCurves} with the classical $\ell^1$-minimization algorithm, then it turns out that the $\ell^1$-minimization curve (for $n=128$, $N=256$) practically coincides with the graph corresponding to the correction of 15 errors ($\approx 12\%$ of measurements located at unknown positions are corrupted ) on Fig. \ref{fig_ErrorCorrectedCurves}.
This observation promotes not only the $\ell^1$-greedy algorithm which is a basic component of Algorithm~B but also supports one of the keystone paradigms of CS stating that the data measurements stored today can be used  much more efficiently in a few years when new algorithms  will be designed.
\par
Even staying within the sparse Gaussian model of errors, we have one parameter whose influence was not considered yet.
This is the magnitude of errors. In classical error correcting codes on the binary field, the magnitude is not an issue. In real-valued encoding, the influence of the magnitude of the error vector on the efficiency of the recovery  is not so obvious. Our experiments with increasing/decreasing the error magnitude in up to 10 times showed that the data reconstruction efficiency increases in all cases even if we do not change the generous weight. There is a common sense explanation of this phenomenon. The Gaussian distribution has the highest information content among distribution with the same variance. The union of the Gaussian data and the Gaussian errors with equal variances has the Gaussian distribution while the weighted union is not the Gaussian one anymore. In fact, it is sparser independently on  whether the weight is less or greater than 1.

%%%%%%%%%%%%%%%%%%%%%%%%%%%%%%
\section{Multisource encoding}
\label{SecMultiSource}
In classical Information Theory the model of the channel with errors is considered  equivalently as a channel with 2 sources of information, i.e., the data transmitted through the channel and undesirable  errors. From this point of view, the errors consume a part of the channel capacity, reducing the capacity of the channel for useful data. This observation leads to the idea to use the error correcting scheme considered in Section \ref{SecVarErrorRate} for encoding/decoding data from a few sources by packing  them in the same output vector. 
\par
To get some idea how efficient this idea can be we consider a model case of 4 sources. In our numerical experiments, each source produces a data vector $\mathbf x_i\in\mathbb R^{64}$. We do not make any special assumption about the sparsity of each specific $\mathbf x_i$. However, we bear in mind that the compound vector $\mathbf x\in\mathbb R^{256}$ has to be somehow sparse to make decoding possible.

\begin{figure}[h]
\centering
\includegraphics[width=3.7in]{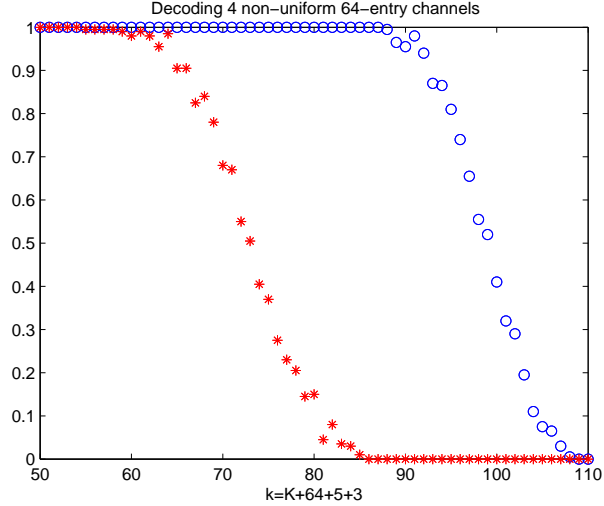}
\caption{Algorithm B with  decoder weights $(3,1,3.6,4)$ vs. Algorithm A for 4 sources with \newline$k_1=K ,\,k_2=64,\,  k_3=5$, and $k_4=3$  entries, 
$\Phi\in\mathbb R^{128\times 256}$ .}
\label{fig_4channel recovery}
\end{figure}

What is especially interesting for us is to outperform Algorithm A when the sources of data contain different amount of information, i.e., when the information is distributed non-uniformly between blocks. In our settings, the input consists of three blocks with the fixed sparsities 64, 5, 3 and one block with the variable sparsity $K$. Thus, the total number of encoded non-zero entries is $k=K+64+5+3$.  For Algorithm B, we set  the block weights equal to $3.0,\,1.0,\,3.6,\,4.0$. We applied only "common sense" in setting the weights. No accurate optimization was performed. The result of the reconstruction is given on Fig.\ref{fig_4channel recovery}.
\par
\par
The numerical experiments show that usage of Algorithm B for decoding data from the multisource input may bring significant benefits over the standard uniform decoding strategy.  The reweighting the blocks completely fits the philosophy of the reweighted optimization. However, the new element in the LGGA consists in combine both the greedy and the generous strategies. While the greedy strategy  sets little penalty weight in the weighted $\ell^1$-norm for the entries suspicious to be used in the encoded data representations, the strategy of Algorithm B quite "generously" changes the weights for entire blocks of information. Setting larger weights in $\ell^1$-norm for the blocks with low information contents, we "temporarily exclude" them from consideration  on initial iterations of Algorithm B, helping blocks with the higher density of the information to be decoded more efficiently.
\par
Now we briefly describe one more possible application, where generous strategy can be useful. In applications above we used partitioning  the matrix $\Phi$ (or $\tF$) generated by a specific problem.  At the same time, in Section \ref{SecVarErrorRate}, we mentioned the phenomenon when
the higher or the lower magnitude of errors increases the capability of LGGA in error correction.
Such phenomenon can be used artificially for increasing the algorithm capability in reconstruction of
sparse data. We give one example showing what benefits are reachable with that approach.
\par
 Let us split the Gaussian matrix $\Phi\in\mathbb R^{n\times N}$ into 2 submatrices $\Phi_1\in \mathbb R^{n\times N_1}$ and $\Phi_2\in\mathbb R^{n\times N_2}$, $N_1+N_2=N$, and 
compose a compound sensing matrix $\Psi=(\Phi_1, \delta \Phi_2)$. For decoding data sensed with the matrix $\Psi$, when $n=N_1=N_2=128$ and $\delta=0.1$, we  introduce little changes in the weight selection of  Algorithm~B. We define the weights as
$$
w_{j,i}:=\left\{\begin{array}{lll}\epsilon, &|x_{j-1,i}|>M_j,&i\le N_1;\\
\epsilon, &|x_{j-1,i}|>M_j/\delta,&i> N_1;\\ 1, &|x_{j-1,i}|\le M_j,&i\le N_1;\\  \lambda, &|x_{j-1,i}|\le M_j,& i>N_1.\end{array}\right.
$$
\par
\begin{figure}[h]
\centering
\includegraphics[width=3.7in]{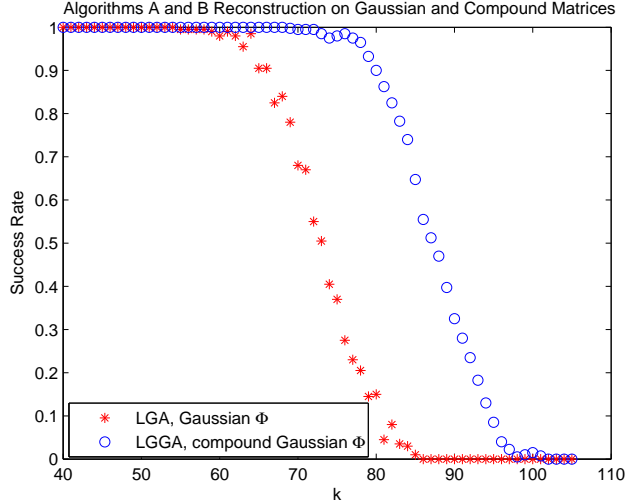}
\caption{Algorithm B on compound matrix vs. Algorithm A on Gaussian matrix, 
$\Phi\in\mathbb R^{128\times 256}$ .}
\label{fig_Couple}
\end{figure}

The result is presented on Fig.  \ref{fig_Couple}.
\par
\par
This experiment obviously shows that the encoding with compound matrices is a very promising tool for 
efficient encoding purposes. At the same time, it should be taken into account that there is a significant stability drawback of that approach.  In Section \ref{SecStability}, we have more extensive discussion about the algorithm stability. However, this is absolutely obvious that the second  half of the data has the theoretical precision of the recnstruction $1/\delta$ times lower than the precision of the first half. So this trick is kind of exchange of   the higher precision for the larger number $k$. From the point of view of settings of formal information theory,  this is not progress at all. Anyway, this idea can be used  when we do not need the same precision for different blocks of information or as a constructive brick for other algorithms, where such compound matrices make a sence (e.g., cf. \cite{KMSSZ}).
%%%%%%%%%%%%%%%%%%%%%%%%%%%%%%%%%%%%%%%%
\section{Adaptive $\ell^1$-Greedy-Generous Decoding Algorithm}
\label{SecAdaptAlg}
In Sections \ref{SecVarErrorRate} and \ref{SecMultiSource} we showed how Algorithm B, applied to linearly compressed data from multiple sources, can efficiently recover the encoded information.
The tuning of Algorithm B uses the side information about the density of the information in each  input sources. 
\par
Embedding this  side information in the code may bring additional (probably non-linear) procedures like protection from the corruption in the channel. While the amount of this side information is tiny, this scheme violates the genuine CS architecture which is linear by nature. Moreover, in most of practical cases, the information density values are not immediately available even for the encoder. Potentially, it may be extracted by the encoder from the raw data at some computational expenses comparable or exceeding CS expenses itself.
Say, to get the sparsity information about an image the wavelet transform has to be applied.  In other cases, like for the channel with errors, the encoder is not aware about the errors in the channel at all. The error rate of a stationary channel can also be measured by sending the empty message  to the decoder.  Thus, to apply Algorithm B we need the distribution of information over the blocks, assumed to be obtained from the encoder,  and information about channel errors, assumed to be measured by the decoder. 
\par
Due to reasonings above, we arrive at conclusion that
the decoding algorithm having an internal estimator of  the amount of information encoded in each block is very desirable. Such estimator has to use only the raw CS encoder output. The block structure, i.e., the partitioning of the vector $\mathbf x$ into the blocks with potentially different information density, also has to be known for the decoder. 
\par
Below, we suggest an adaptive algorithm working on two sources of information. 
While it is quite efficient within the considered settings, it serves as just  one successful example on the way to a really universal adaptive algorithm acting on the variety of block configurations as well as on the different information contents of the blocks.
\par
We conducted some research bringing the algorithm for adapting weights between blocks, according to 
intermediate internal estimates made on LGGA iterations.  Due to iterative nature of LGA we do not try to find the optimal weights for the blocks from the beginning. We rather dynamically update those weights in the direction of the higher sparsity of the result. 
For our estimates, we  compute  0.5-quasi-norm $\|\mathbf x\|_{0.5}:=(\sum \sqrt{|x_i|})^2$ normalized by the Euclidean norm $\|\mathbf x\|_2:=\sqrt{\sum x_i^2}$:
 $$
s(\mathbf x):=\|\mathbf x\|_{0.5}/\|\mathbf x\|_2.
$$
Since the vector $\mathbf x$ is unknown in the beginning of the decoding, we will use heuristic methods measuring the relative  behavior of this characteristic for the blocks of the vector $\mathbf x$ on consecutive iterations of LGGA. 
\par
We consider only the simplest case of two  blocks of information with potentially different densities. We assume that the blocks have equal size, i.e., for modeling settings 
from previous section when $N=256$, we assume that the block sizes $N_1$ and $N_2$ are equal to 128. Let us agree that the index sets are  $\{1,\dots,N/2\}$ and $\{N/2+1,\dots, N\}$ for the first and for the second blocks correspondingly. We will use the superscripts 1 and 2 for denoting the corresponding  subvectors and the subscripts $"p"$ and $"c"$ for the minimum of $\ell^2$-norm  (pseudoinverse) solution to $\Phi \mathbf x=\mathbf y$ and the current estimate of $\mathbf x$ made by LGGA. 
\par
The next iteration weights  will be functions of the value
$$
S=\frac{s(\mathbf x^1_c)\cdot s(\mathbf x^2_p)}{s(\mathbf x^2_c)\cdot s(\mathbf x^1_p)}.
$$
Our heuristic approach is based on the experimental observation  that when the second block has lower
density of the information and the algorithm transforming $\mathbf x_p$ into $\mathbf x_c$ has some sparsifying properties, then $S>1$ at least in the sense of the expectation. The value of 
$S$ grows with the growth of the non-uniformity of the information distribution between the blocks.
\par
Accepting those results, we have to introduce a function  setting the block weights in the greedy-generous algorithm.  We set the weight $W_1=1$ for the first block  and find the second weight as 
\begin{equation}
\label{WEIGHT}
W_2(S)=\left\{\begin{array}{ll}
0.65\sqrt{S}+0.35\sqrt[4]{S}, &S\ge 1;\\
(0.65\sqrt{S}+0.35\sqrt[4]{S})^{-1}, & S<1.
\end{array}
\right.
\end{equation}
\par

{\bf Algorithm C}
({\it adaptive $\ell^1$-Greedy-Generous Algorithm})
\par
\begin{enumerate}
\item{\it Preliminaries}
\begin{enumerate}
\item
Compute the minimum $\ell^2$-norm solution $\mathbf x_p:=\Phi^T(\Phi\Phi^T)^{-1}\mathbf y$.
\item
Compute $s(\mathbf x_p^1)$ and $s(\mathbf x_p^2)$
\end{enumerate}
\item {\it Initialization}.
\begin{enumerate}
\item
Set $w_{0,i}\equiv 1$, $i=1,\dots,N$, $W_1=1$.
\item
 Find the solution  $\mathbf x_0$ to (\ref{sle}) providing the minimum for (\ref{weighted}).
\item
 Set $M_0:=\max\{|x_{0,i}|\}$ and $j:=1$.
\end{enumerate}

\item
{\it General ALGGA Step}.
\begin{enumerate}
\item
Find $S$ and $W_2=W_2(S)$
\item Set $M_j:=\alpha M_{j-1}$ and the weights
$$
w_{j,i}:=\left\{\begin{array}{ll}\epsilon, &|x_{j-1,i}|>M_j,\\ W_1, &|x_{j-1,i}|\le M_j,\,i\le N/2 \\W_2, &|x_{j-1,i}|\le M_j,\,i>N/2.\end{array}\right.
$$
\item
 Find the solution  $\mathbf x_j$ to (\ref{sle}) providing the minimum for (\ref{weighted}).
\item
 Set $j:=j+1$.
\item
If Stopping Criterion is not satisfied, GOTO 3 (a).
\end{enumerate}
\end{enumerate}
\par

\begin{figure}[h]
\centering
\includegraphics[width=3.7in]{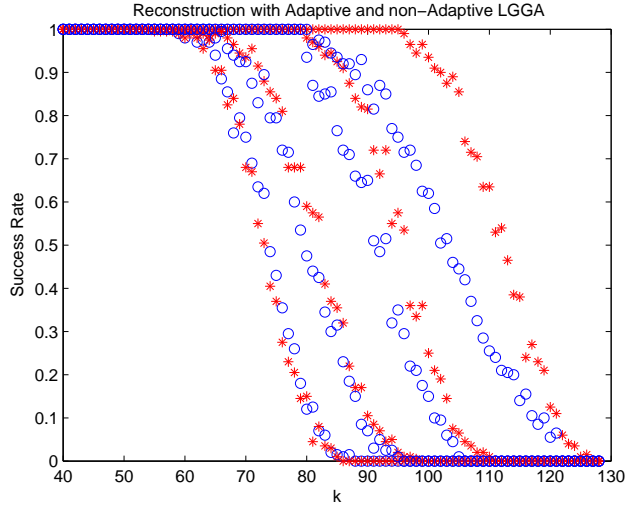}
\caption{Algorithm C ('o', adaptive) vs. Algorithm B ('*', optimal) success rates\newline  for 2-channel input $k=k_1+k_2$,\, $k_2= 37,\, 15,\,5,\,1$  (from left to right).}
\label{fig_5AdaptiveNonAdaptive}
\end{figure}

The results of our numerical experiment are presented on Fig. \ref{fig_5AdaptiveNonAdaptive}. For all graphs the horizontal coordinate $k$ corresponds to the total number of non-zero entries in both blocks. The graphs plotted 
with "*" correspond to the recovery with Algorithm B when the density of the information is known
before decoding, whereas the graphs plotted with "$\circ$" reflects the results obtained with adaptive Algorithm~C. There are four pairs of the graphs. 

\par

The most left pair corresponds to the uniform information distribution between blocks. 

To be more precise, we set a number of non-zero entries in block 2 as 37 and a number $k_1$ of non-zero entries in block 1 is changing, giving  the total number $k=k_1+37$. The number 37 is one half of the value 
corresponding to the level 0.5 of the success curve for Algoritm A. Thus, in the neighborhood of
$k=74=2\times 37$ we have non-zero entries uniformly distributed over the entire range $1\div 256$. Of course, in this case,  the optimal settings  for Algorithm B are $W_1=W_2=1$. Therefore Algorithms A and B have identical output.
\par
Since those 2 most left curves are practically coincide, this means that Algorithm C makes a reasonable decision, setting $W_2\approx 1$ and showing its stability.
\par
The second pair of  curves from the right corresponds the case $k_2=15$. The Algorithm C output is very close to the output of Algorithm B (with $W_2=1.7$).  For the case $k_2=5$ (the third pair of curves), Algorithm C gives very big advantage over Algorithm A (the case of uniform distribution), however it concedes visually to Algorithm B (with $W_2=2.5$). For the case $k_2=1$ (the most right pair of graphs), Algorithm C is significantly less efficient than Algorithm B (with $W_2=6.5$) but has a huge advantage over Algorithm A.  
\par
We would like to emphasize that the function defining the weight $W_2$ in (\ref{WEIGHT}) does not pretend  to be optimal or somehow universal. We just wanted to show that the idea of efficient decoding of multichannel code with unknown characteristics of channels may have a quite satisfactory adaptive solution.
Anyway, we believe that the significant part of Algorithm C losses for highly non-uniform distributions is mainly because of statistical inconsistency of the estimate of S rather than imperfection of formula
(\ref{WEIGHT}). This is obvious that a few non-zero values have insignificant influence on the values
of sparsity estimates, especially for pseudoinverse solutions $s(\mathbf x_p^j)$. Therefore the significant fluctuations of $S$ may mislead formula (\ref{WEIGHT}) in setting $W_2$.
\par
Experiments with change of the variance of $\mathbf x^2$ give an argument in favor of our claim above. Indeed, we introduce a factors for $\mathbf x^2$ in the range $0.1\div10$ and checked the reconstruction with Algoritm C. In all cases, the higher rate of recovery was observed.  In particular, for the factor 10, the curve for $k_2=1$  became practically coinciding with the Algorithm B curve from 
Fig. \ref{fig_5AdaptiveNonAdaptive}. We suppose that the reason of this improvement is the increased consistency of the sparsity estimate.
\par
\section{$\ell^1$-Greedy-Generous Algorithm Stability.}
\label{SecStability}
The stability of the algorithm with respect to the noise in measurements is a necessary requirement 
for its applicability to the real world problems. The noise may have very different reasons. The most typical forms of  noise are noise of transmission and noise of quantization.
Any noise  reduces the precision of the result. However, technically, the goal of this section is to estimate the stability of the fact of recovery of sparse  representations with reasonable precision rather then 
fighting for the best possible precision itself. The CS decoding is  a highly non-linear operation. So  the stability means for CS
much more than just the rate of the dependence of $\|\mathbf x-\hat{\mathbf x}\|$ from
 $\|\mathbf y-\hat{\mathbf y}\|$ because even the continuity of such dependence is under  question. Generally, no algorithm can be stable for all $k$ not exceeding $n-1$.
Indeed, each value of the measurement precision has its own theoretical maximum level sparsity $k$  admitting the recovery. 
In particular, the case $k=n-1$ requires infinitely high precision when $N/n\ge\alpha>1$; $n\to\infty$. For $k/n=\beta<1$, the reconstruction with moderate precision is possible up to some
 level of noise in $\mathbf y$. For the noise/precision level above that threshold, no precision of reconstruction can be guaranteed.  The closer $\beta$ to 1, the lower the threshold is.
We checked the result of reconstruction for many different levels of noise. 
The greedy-generous algorithm has shown the high stability to the noise in the measurements 
$\mathbf y$. Let us describe the settings of numeric experiments supporting the stability claim.
\par
In the noisy environment, we run the same adaptive LGGA with the same parameters but we change the criterion 
of success. The threshold value  for "success" was set as $\delta=1.5\cdot2^{NH(k/N)/n}\sqrt{n}\sigma$, where 
$$
H(p):=-p\log_2p-(1-p)\log_2(1-p),
$$
$\sigma$ is the noice variance. 
When $\|\hat{\mathbf x}-\mathbf x\|<\delta$, where $\hat{\mathbf x}$ is the estimate obtained with the algorithm, we say that the reconstruction of ${\mathbf x}$ is successful. 
\par
We use an indirect approach for setting the threshold $\delta$. The reasoning is as follows.
Due to normalization of the matrix $\Phi$, its random $n\times n$ submatrix is "almost" unitary.
So in the best case scenario any deviation from the correct values in the $\mathbf y$-domain will be
conversed into the same deviation in  $\mathbf y$-domain. Let us obtain a numerical estimate for this conversion.
First of all, the noise level $\sigma$ at $n$ entries gives the total noise $\sqrt{n}\sigma$ at  $\mathbf y$.
The level of the noise can be translated into a number of significant bits in the entries of $\mathbf y$.
Unfortunately, we cannot hope that all those bits to be used for representation of components of the vector $\mathbf x$. Indeed, the linear encoding formula (\ref{sle}) assumes that the data vector 
$\mathbf x$ can be recovered from the measurements $\mathbf y$. In particular, this means that 
the vector $\mathbf y$ has to contain not only information about bits in digital representation of 
$\mathbf x$ but also (maybe implicitely) information about indices of non-zero entries.  For encoding that information
we need $H(k/N)$ bits per entry. The total number of entries in $\mathbf x$ is $N$. The information about indices has to be put into the vector $\mathbf y$ consisting of $n$ entries, i.e., each entry of $\mathbf y$ has to "reserve"  $NH(k/N)/n$ bits for encoding the indices. Those bits iminently decrease the precision of the reconstruction at least in  $2^{NH(k/N)/n}$. Thus, we cannot expect  the precision of reconstruction  $\|\hat{\mathbf x}-\mathbf x\|$  higher than $2^{NH(k/N)/n}\sqrt{n}\sigma$. Numerical experiments shows that this value is really very close to the precision of the estimate $\hat{\mathbf x}$. We introduce the additional factor 1.5 just to avoid unfair "losses" of the success when the precision is slightly higher due to some random fluctuations.
\par
\begin{figure}[h]
\centering
\includegraphics[width=3.7in]{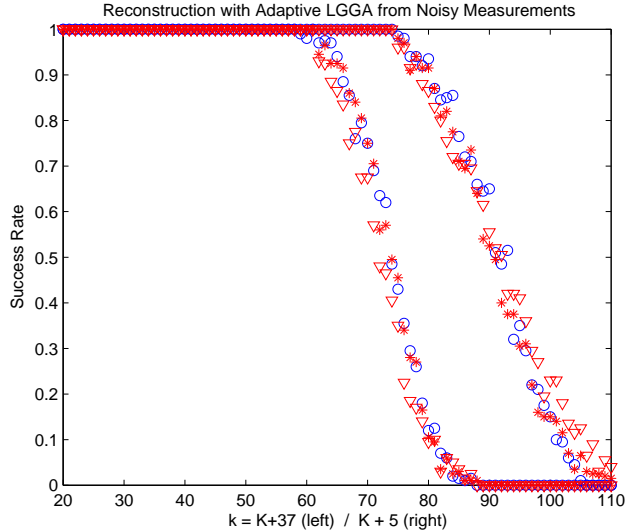}
\caption{Algorithm C success rate  for noise $\sigma=0\, (\circ), \,0.01\, (*), 0.03 \,(\nabla)$, \newline 
for $k_2=37$ (left) and $k_2=5$ (right).}
\label{fig_5NoisyMeasurments}
\end{figure}
\par
The graphs of the success rate  for non-uniform distribution of the information with a fixed number
$k_2=5$ (right triplet) and $k_2=37$ (right triplet) for $\sigma=0,\, 0.01,\, 0.03$ are given on Fig.~\ref{fig_5NoisyMeasurments}. 
\par
The exciting fact that "the best case scenario" estimate turns out to be a good tool for the prediction 
of the deviation of the algorithm output from the true vector is very  optimistic for CS perspectives.
Indeed, the threshold estimate is based on elementary information theory principles which cannot 
be overcome with any algorithm.  Compression of information allowing this precision of the reconstruction can be achieved if we separately encode the  digital (quantized) information about
the entries and the information about indices. The index part of the information can be encoded with the optimal bit budget by
the arithmetic encoder whose output is extremely vulnerable to the errors. In fact, the results of numerical modeling shows that linear compressor (\ref{sle}) encodes
the combination  of digital bits and indices in very wise form. The indices are encoded implicitly.
They are a part of Joint Source-Channel linear encoding.  Within a wide range of losses of precision in the entry representation, reconstruction does not destroy the structure (index information), providing 
decoding with the quality progressively depending on the channel noise.

%%%%%%%%%%%%%%%%%%%%%%%%%%%%%%%%%%%%%%%%%%
\section{Conclusions}
It is shown that Compressed Sensing encoding has very strong internal error correcting capability. No special redundant encoding is required for further error correction. Essentially, the error protection capability is defined by the amount of the space underloaded with data by the encoder. The less the  data information content, the higher natural error protection level of the code. Thus,   Compressed Sensing provides highly efficient Joint Source-Channel Encoding method. 
\par
We showed that a minor modification of the $\ell^1$-greedy decoding algorithm, which we call the  $\ell^1$-greedy-generous algorithm (LGGA), allows to correct errors efficiently.
\par
Using the well-known parallel between combination of data and errors in the channel with a few sources joint encoding, we applied the same decoding algorithm to  decoding the multi-source code. 
In the case of non-uniform informational contribution of the sources into the encoder output, our new $\ell^1$-greedy-generous algorithm significantly outperforms all known algorithms,  including $\ell^1$-greedy, in recovering Gaussian data encoded with Gaussian matrices.
\par
While the knowledge of approximate distribution of information between blocks is very desirable, this distribution can be estimated dynamically during iterations of LGGA. We suggest algorithm automatically adapting the "generous" weights of the $\ell^1$-greedy-generous algorithm to input with unknown (possibly non-uniform) distribution of information between blocks.

\end{document}